\documentclass{article}

\usepackage{amsmath, amssymb, amsfonts, amsthm, bm, eucal, graphicx, overpic}

\DeclareMathOperator{\Area}{Area}
\DeclareMathOperator{\sdet}{Sdet}
\newcommand{\epsi}{\varepsilon}
\newcommand{\norm}[1]{\|#1\|}
\newcommand{\Norm}[1]{\bigl\|#1\bigr\|}
\newcommand{\ipr}[2]{\langle #1, #2 \rangle}
\newcommand{\Tr}{\mathsf{T}}
\newcommand{\Quat}{\mathbb{H}}
\newcommand{\Q}{\mathbb{Q}}
\newcommand{\R}{\mathbb{R}}
\newcommand{\Z}{\mathbb{Z}}
\newcommand{\Proj}{\mathbb{P}}
\newcommand{\C}{\mathbb{C}}
\newcommand{\D}{\mathbb{D}}
\newcommand{\vect}[1]{\bm{#1}}
\newcommand{\va}{\vect{a}}
\newcommand{\vb}{\vect{b}}
\newcommand{\vc}{\vect{c}}

\newcommand{\vo}{\vect{o}}

\newcommand{\vv}{\vect{v}}

\newcommand{\collection}[1]{{\mathcal#1}}
\newcommand{\dualS}{\collection{S}}
\newcommand{\CS}{\collection{S}}
\newcommand{\abs}[1]{|#1|}
  \newcommand{\Abs}[1]{\left|#1\right|}
\newcommand{\myangle}{\sphericalangle}
  \def\0{^{\phantom0}}
  \def\9{_{\phantom9}}

\theoremstyle{plain}
\newtheorem{theorem}{Theorem}
\newtheorem{lemma}{Lemma}

\begin{document}

\title{Sylvester-Gallai Theorems for Complex Numbers and Quaternions}
\author{Noam Elkies\thanks{Department of Mathematics, Harvard University, 
        Cambridge MA 02138.
        E-mail: \texttt{<elkies@math.harvard.edu>}.
	Supported in part by NSF grant DMS-0200687.}
\and
        Lou M. Pretorius\thanks{Department of Mathematics and Applied Mathematics,
        University of Pretoria,
        Pretoria 0002, South Africa.
        E-mail: \texttt{<lpretor@scientia.up.ac.za>}.
	Supported in part by the National Research Foundation 
	under Grant number 2053752.}
\and
        Konrad J. Swanepoel\thanks{Department of Mathematical Sciences,
	University of South Africa, PO Box 392,
	Pretoria 0003, South Africa.
	E-mail: \texttt{<swanekj@unisa.ac.za>}.
        Supported by the National Research Foundation 
	under Grant number 2053752.}}
\date{}
\maketitle

\begin{abstract}
A Sylvester-Gallai (SG) configuration is a finite set~$S$\/ of points
such that the line through any two points in $S$\/
contains a third point of $S$.
According to the Sylvester-Gallai Theorem,
an SG configuration in real projective space must be collinear.
A problem of Serre (1966) asks whether an SG configuration
in a complex projective space must be coplanar.
This was proved by Kelly (1986) using a deep inequality of Hirzebruch.
We give an elementary proof of this result, and then extend it to show
that an SG configuration in projective space over the quaternions
must be contained in a three-dimensional flat.
\end{abstract}

\section{Introduction}
We denote the fields of real and complex numbers by $\R$ and $\C$,
respectively, and the division ring of quaternions by $\Quat$.
We let $\Proj^n(\D)$ denote the $n$-dimensional projective space
over the division ring $\D$.
A finite subset~$S$\/ of $\Proj^n(\D)$
is a \emph{Sylvester-Gallai configuration} (SG configuration)
if for any distinct $x,y\in S$\/ there exists $z\in S$\/
such that $x,y,z$ are distinct collinear points.
It is a classical fact that the nine inflection points
of a non-degenerate cubic curve in $\Proj^2(\C)$
constitute an SG configuration.  These nine points are not collinear.
Sylvester \cite{Sylvester} asked whether
an SG configuration in $\Proj^2(\R)$ must be collinear.
Later Erd\H{o}s independently asked this question \cite{MR84e:01073},
which was solved by Gallai and others \cite{Steinberg}
(see \cite{MR92b:52010} for a survey).
Since the $n$-dimensional case trivially follows from
the two-dimensional case, we formulate this result as follows:

\begin{theorem}[Sylvester-Gallai]\label{sylvester-gallai}
Every SG configuration in $\Proj^n(\R)$ is collinear.
\end{theorem}

Serre \cite{Serre} asked whether an SG configuration in $\Proj^n(\C)$
must be coplanar (i.e., must lie in a two-dimensional flat).
This was solved by Kelly \cite{MR87k:14047}
using an inequality of Hirzebruch \cite{MR84m:14037}
involving the number of incidences of points and lines in $\Proj^2(\C)$.
This inequality follows from deep results in algebraic geometry.

\begin{theorem}[Kelly]\label{kelly-elkies}
Every SG configuration in $\Proj^n(\C)$ is coplanar.
\end{theorem}

We show how one of the many proofs of Theorem~\ref{sylvester-gallai}
may be generalized to give an elementary proof of the above theorem.

Of course a similar question may be asked for projective space
over the quaternions: What is the smallest dimension~$k$\/
such that every SG configuration in $\Proj^n(\Quat)$
must lie in a $k$-dimensional flat?
We show that the elementary proof of Theorem~\ref{kelly-elkies}
also generalizes to this case, giving the following:
\begin{theorem}\label{main}
Every SG configuration in $\Proj^n(\Quat)$
lies in a three-dimensional flat.
\end{theorem}
We do not know whether this result is sharp:
we have no example of an SG configuration that spans $\Proj^3(\Quat)$.

In Sections~\ref{complexsection} and \ref{quaternionsection}
we prove Theorems~\ref{kelly-elkies} and \ref{main}, respectively.
We discuss the proof of Theorem~\ref{sylvester-gallai} in detail
in Section~\ref{realsection}, since it forms the basis of the proofs
of Theorems~\ref{kelly-elkies} and \ref{main}.
In the next section we discuss the mathematics
needed for the above proofs.

\section{Notation and definitions}
We consider the $(k+1)$-dimensional subspaces of $\D^{n+1}$
to be the $k$-dimensional flats of the $n$-dimensional projective space
$\Proj^n(\D)$ over a division ring $\D$ in the usual way.
When passing to the affine space $\D^n$ we use barycentric coordinates,
i.e., we take the hyperplane at infinity of $\Proj^n(\D)$
to correspond to $\sum_{p=1}^{n+1}x_p = 0$ in $\D^{n+1}$,
and the coordinates of a point $P$\/ not at infinity
to be the coordinates of the intersection of the hyperplane
$\sum_{p=1}^{n+1}x_p = 1$
with the $1$-dimensional subspace of $\D^{n+1}$ corresponding to~$P$.
To each $\Proj^n(\D)$ we may associate its dual $\Proj^n(\D)^\ast$,
with a duality map taking $k$-flats in $\Proj^n(\D)$ to $(n-1-k)$-flats
in $\Proj^n(\D)^\ast$.  In particular, the dual of a point
is an $(n-1)$-flat or \emph{hyperplane}.
When constructing the dual we let $\D^{n+1}$ be the right vector space
of column vectors over~$\D$.
Then $\Proj^n(\D)^\ast$ is the projective space
coming from the dual $(\D^{n+1})^\ast$,
which is the left vector space of row vectors over $\D$.
See \cite{Samuel} for more on projective spaces.

By duality we may associate with any SG configuration~$S$\/
its \emph{dual SG configuration} $\dualS$ in the dual projective space.
A dual SG configuration in $\Proj^n(\D)$ is then a finite set of hyperplanes
such that for any distinct $\Pi_1,\Pi_2\in\dualS$
there exists $\Pi_3\in\dualS$ such that
$\Pi_1,\Pi_2,\Pi_3$ are distinct and
$\Pi_1\cap\Pi_2\cap\Pi_3$ is an $(n-2)$-flat.

We use the usual representation $\alpha=t+xi+yj+zk=t+\vv$
for a quaternion, with $\vv$ a vector in $\R^3$.
We let $\abs{\alpha}$ denote the norm of~$\alpha$,
and $\alpha^\ast=t-\vv$ its conjugate.
We consider the $n$-dimensional vector space $\Quat^n$
to be the space of column vectors,
with scalar multiplication from the right,
and thus the action of linear transformations
as matrix multiplication from the left.

It is well-known that $\Quat$ can be represented
in the ring of $2\times2$ complex matrices,
by identifying $a+bi+cj+dk = (a+bi) + (c+di)j$ with
\[ \left[
\begin{array}{cc}
a+bi & -c+di\\
c+di & \;\,\,a-bi
\end{array}\right].\]
If we replace each entry $a_{pq}$ of an $n\times n$ matrix~$A$\/
with quaternion entries by its corresponding $2\times 2$ complex matrix,
we obtain a $2n\times 2n$ complex matrix $A_\C$.
The \emph{Study determinant} of $A$ is then defined by
\[ \sdet(A) = \det(A_\C).\]
See \cite{MR97j:16028} for an exposition.
The Study determinant is a non-negative real number.
It is multiplicative: $\sdet(AB)=\sdet(A)\sdet(B)$.
Hence $\sdet(A)>0$ if $A$ is invertible.
The Study determinant of the $1\times1$ matrix $[\alpha]$
is $\sdet[\alpha]=\abs{\alpha}^2$.

We let $\ipr{\cdot}{\cdot}$ denote the standard inner product in $\R^n$.

\section{Proof of Theorem~\ref{sylvester-gallai}}\label{realsection}
We review a known\footnote{
  Known, but hard to locate in the literature\ldots
  \ One of us (Elkies) recalls reading, around 1980,
  a version of this proof in which Lemma~\ref{trianglelemma}
  is proved by reducing to the case that $\Delta$ is equilateral
  and using the largest angle of the circumscribing triangle
  to obtain another triangle of equal or smaller area, as in \cite{AMMProblem}.
  But we cannot find a reference for this argument.
  Several proofs use the dual line configuration $\dualS$,
  but then conclude by applying either Euler's formula
  or the order properties of~$\R$, neither of which
  can be used over~$\C$ or~$\Quat$.
  }
proof of the Sylvester-Gallai theorem (Theorem~\ref{sylvester-gallai}).
It is sufficient to prove the case $n=2$,
since if there is an SG configuration $S$\/
spanning a flat of dimension more than two, we may choose three points
$a,b,c\in S$ spanning a \hbox{$2$-flat} $\Pi$, and then
$S\cap\Pi$ is again an SG configuration in the projective plane $\Pi$.
Thus we assume that we have an SG configuration $S$ spanning $\Proj^2(\R)$.

Let $\dualS$ be its dual SG configuration.
Then $\bigcap_{l\in\dualS} l = \emptyset$.
Choose a line at infinity
that is not one of the lines in $\dualS$ and does not pass
through any point of intersection of two lines in $\dualS$.
Then in the affine plane obtained by removing the line at infinity,
$\dualS$ is a dual SG configuration of mutually non-parallel lines.

\begin{figure}
\begin{center}
\begin{overpic}[scale=0.6]{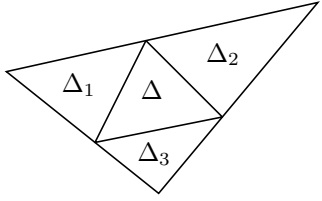}
\put(20,33){$\Delta_1$}
\put(62,42){$\Delta_2$}
\put(42,13){$\Delta_3$}
\put(43,31){$\Delta$}
\end{overpic}
\end{center}
\caption{$\Area(\Delta_i)\leq\Area(\Delta)$ for some $i$}\label{fig1}
\end{figure}

Since $\dualS$ does not have a common point, $\dualS$ contains
some three non-concurrent lines determining a triangle.
Among all such triangles, choose one of minimum area.
Let the lines of this triangle be $\ell_1,\ell_2,\ell_3$.
For any $i,j\in\{1,2,3\}$ with $i<j$, there is a third line
$\ell_{ij}\in\dualS$ passing through the intersection point
of $\ell_i$ and $\ell_j$.
Then we obtain a contradiction from the following result,
which is a reformulation of a well-known geometric inequality 
first published by Debrunner \cite{Debrunner}; see also \cite{AMMProblem}.
See \cite[inequalities~9.1, 9.2 and 9.3]{Bottema} for a discussion 
of this and related inequalities.
\begin{lemma}\label{trianglelemma}
One of the nine triangles
\[ \begin{array}{lll} 
\ell_{12}\ell_2\ell_3 & \ell_1\ell_{12}\ell_3 & \ell_1\ell_2\ell_{13} \\
\ell_{13}\ell_2\ell_3 & \ell_1\ell_{23}\ell_3 & \ell_1\ell_2\ell_{23} \\
\ell_1\ell_{12}\ell_{13} & \ell_{12}\ell_2\ell_{23} & \ell_{13}\ell_{23}\ell_3
\end{array}\]
has area at most that of triangle $\ell_1\ell_2\ell_3$.
Furthermore, if none of them has area
strictly less than that of triangle $\ell_1\ell_2\ell_3$,
then we must have parallel lines
$\ell_1\parallel\ell_{23}$,
$\ell_2\parallel\ell_{13}$, and
$\ell_3\parallel\ell_{12}$.
\end{lemma}
See Figure~\ref{fig1}.
In order to make the next two proofs easier to read,
we now prove this lemma in a complicated way.

\begin{proof}[Proof of Lemma~\ref{trianglelemma}]
We identify the affine plane with the plane $x_1+x_2+x_3=1$ in $\R^3$.
Then $\ell_i$ has equation $x_i=0$, and
the vertices of the triangle $\ell_1\ell_2\ell_3$ are
$P_1=[1,0,0]^{\Tr}$, $P_2=[0,1,0]^{\Tr}$, $P_3=[0,0,1]^{\Tr}$.
Since $P_k\in\ell_{ij}$ when $k\neq i$ and $k\neq j$,
the line $\ell_{ij}$ has an equation of the form
$\alpha x_i+\beta x_j=0$
for some $\alpha,\beta$, with $\alpha,\beta\neq 0$
since $\ell_{ij}\neq\ell_i,\ell_j$.
We now set $\ell_{ji}=\ell_{ij}$.
Then for any distinct $i,j$, we can write the equation of $\ell_{ij}$
as $\alpha_{ij}x_i+x_j=0$, where 
\begin{equation}
\label{realone}
\alpha_{ij}\alpha_{ji}=1\text{ for all distinct $i,j \in \{1,2,3\}$.}
\end{equation}
Note that for each permutation $(i,j,k)$ of $\{1,2,3\}$ we have
$\alpha_{ij}=1$ if and only if $\ell_{ij}\parallel\ell_k$.

We now consider the areas of the different triangles.
It is easily seen that the area of a triangle with vertices $P,Q,R$\/
is $\frac{\sqrt{3}}{2}\Abs{\det[P,Q,R]}$.
Using the equations of the lines of triangle $\ell_1\ell_{12}\ell_{3}$,
we compute that its vertices
$P_1'=\ell_{12}\cap\ell_3,
P_2'=\ell_1\cap\ell_3,
P_3'=\ell_1\cap\ell_{12}$
have coordinates
\[
P_1'=[(1-\alpha_{12})^{-1}, -\alpha_{12}(1-\alpha_{12})^{-1},0]^{\Tr},
\
P_2'=P_2\0,
\
P_3'=P_3\0,
\]
whence $\Abs{\det[P_1',P_2',P_3']}=\abs{1-\alpha_{12}}^{-1}$.
Since $\det[P_1,P_2,P_3]=1$, the assumption that $\ell_1 \ell_2 \ell_3$
has area no larger than $\ell_1 \ell_{12} \ell_3$ implies
$\abs{(1-\alpha_{12})^{-1}}\geq 1$, that is,
\[
\abs{1-\alpha_{12}}\leq 1.
\]
Likewise, by calculating the area of $\ell_1\ell_{12}\ell_{13}$
we find that this triangle has area less
than or equal to the area of $\ell_1 \ell_2 \ell_3$
if and only if
\[
\abs{1-\alpha_{12}-\alpha_{13}}\leq 1.
\]
By permuting indices we obtain the following $9$ inequalities
in the $6$ real variables $\alpha_{ij}$:
\begin{align}
\label{realtwo}
\abs{1-\alpha_{ij}}\leq 1,& \qquad i,j\text{ distinct};\\
\label{realthree}
\abs{1-\alpha_{ij}-\alpha_{ik}}\leq 1,&\qquad i,j,k\text{ distinct.}
\end{align}
By \eqref{realtwo}, all $\alpha_{ij}\geq 0$.
By \eqref{realone} and the AGM inequality,
$\alpha_{ij}+\alpha_{ji}\geq 2$,
with equality if and only if $\alpha_{ij}=\alpha_{ji}=1$.
Summing over the three pairs $\{i,j\}$, we deduce
\[
\sum_{i,j: i\neq j}\alpha_{ij}\geq 6.
\]
On the other hand, from \eqref{realthree} we obtain
$\alpha_{ij}+\alpha_{ik}\leq 2$.
Summing over all $i$ we obtain
\[
\sum_{i,j: i\neq j}\alpha_{ij}\leq 6.
\]
Therefore, $\sum_{i\neq j}\alpha_{ij}=6$,
and we conclude that all $\alpha_{ij}=1$,
implying that $\ell_{ij}\parallel\ell_k$
for each permutation $(i,j,k)$ of $\{1,2,3\}$.
\end{proof}

This finishes the proof of Theorem~\ref{sylvester-gallai}. \qed

\begin{figure}
\begin{center}
\includegraphics[scale=0.5]{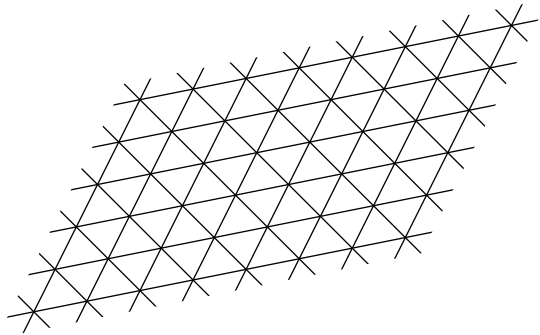}
\end{center}
\caption{Skew triangular lattice}\label{trilattice}
\end{figure}

\bigskip
A slightly different and harder, but more suggestive, proof
could be obtained even if we did not choose the line at infinity
so that no two lines of~$\dualS$ are parallel.
Then Lemma~\ref{trianglelemma} would not immediately yield
a contradiction, but we could apply it again to each of the triangles
$\Delta_i$ of Figure~\ref{fig1} to obtain yet more parallel lines
in~$\dualS$.  Proceeding inductively, we would find that
$\dualS$ contains the line $x_i=m$ for each $i\in\{1,2,3\}$ and $m\in\Z$.
See Figure~\ref{trilattice}.
This is a contradiction because $\dualS$ must be finite --- though
the infinite family of lines $x_i=m$ is locally finite
and does satisfy the dual SG incidence condition:
no point in the plane is contained in exactly two of these lines.
Up to affine linear transformation, this is thus the unique
such configuration with a triangle of minimal area.

There are two known relatives of Lemma~\ref{trianglelemma} 
\cite[inequalities~9.2 and 9.3]{Bottema}, one where area is replaced by 
perimeter,
and the other by the inradius of the triangle.
These two related inequalities may be used in the same way as 
Lemma~\ref{trianglelemma} to deduce Theorem~\ref{sylvester-gallai}.

\section{Proof of Theorem~\ref{kelly-elkies}}\label{complexsection}
We argue as in the proof of Theorem~\ref{sylvester-gallai}.
It is again clearly sufficient to prove the case $n=3$.
As before, we dualize.
Let $\dualS$, then, be a dual SG configuration in $\Proj^3(\C)$
with empty intersection, and as before we choose the plane at infinity
to be one not in $\dualS$ and not containing any line of intersection
of two planes in $\dualS$.
Then in the resulting affine space,
which we identify with the hyperplane $x_1+x_2+x_3+x_4=1$ in $\C^4$,
the configuration $\dualS$ consists of pairwise non-parallel planes.
By hypothesis, $\dualS$ contains at least one four-tuple of planes
$\Pi_1,\Pi_2,\Pi_3,\Pi_4$ with empty intersection.
We call any four such planes a \emph{tetrahedron}
with \emph{vertices} $P_j=\bigcap_{k\neq j}\Pi_k$, $j=1,2,3,4$.
As before, we may compare ``volumes'' of tetrahedra, where
we now define the \emph{volume} of a tetrahedron $\Pi_1\Pi_2\Pi_3\Pi_4$
to be $\Abs{\det[P_1, P_2, P_3, P_4]}$.
(The ``volume'' could be defined to be some constant positive multiple
of this determinant, without changing the argument, as the factor
$\sqrt{3}/2$ had no effect in our proof of Lemma~\ref{trianglelemma}
in $\Proj^2(\R)$.)
Fix a tetrahedron $\Pi_1\Pi_2\Pi_3\Pi_4$ of minimum volume,
and choose coordinates $x_1,x_2,x_3,x_4$
so that the vertices $P_j$ of this tetrahedron
are the standard unit basis vectors
$[\delta_{j1},\delta_{j2},\delta_{j3},\delta_{j4}]^{\Tr}$ of~$\C^4$.
Now $\Pi_j$ has equation $x_j=0$.
By the multiplicativity of determinants,
$\Pi_1\Pi_2\Pi_3\Pi_4$ remains a tetrahedron of minimum volume,
which is now~$1$.

For any $j,k\in\{1,2,3,4\}$ with $j<k$, choose a plane
$\Pi_{jk}$ in $\dualS$ with $\Pi_j\cap\Pi_k\subset\Pi_{jk}$.
Then $\Pi_{jk}$ has an equation
$\alpha x_j+\beta x_k=0$ with $\alpha,\beta\neq 0$, since
$P_m\in\Pi_{jk}$ for all $m\neq j,k$ and $\Pi_{jk}\neq\Pi_j,\Pi_k$.
Thus we may write this equation as $x_j+\alpha_{jk}x_k=0$,
and the equation of $\Pi_{kj}:=\Pi_{jk}$ as $x_k+\alpha_{kj}x_j=0$,
where
\begin{equation}
\label{complexone}
\alpha_{jk}\alpha_{kj}=1\text{ for all distinct $j,k\in\{1,2,3,4\}$.}
\end{equation}
We now compare the volume of $\Pi_1\Pi_2\Pi_3\Pi_4$
to the volumes of other tetrahedra.
For example, $\Pi_1\Pi_{12}\Pi_3\Pi_4$ has vertices
$P'_1, P'_2, P'_3, P'_4$, where
\[
P_1'=[(1-\alpha_{12})^{-1},-\alpha_{12}(1-\alpha_{12})^{-1},0,0]^{\Tr},
\
P'_2 = P_2\0,
\
P'_3 = P_3\0,
\
P'_4 = P_4\0.
\]
Then $\Abs{\det[P_1',P_2',P_3',P_4']}=\abs{(1-\alpha_{12})^{-1}}$,
and we obtain as before
\[
\abs{1-\alpha_{12}}\leq 1.
\]
By considering $\Pi_1\Pi_{12}\Pi_{13}\Pi_4$ we obtain
\[ \abs{1-\alpha_{12}-\alpha_{13}}\leq 1,\]
and by considering $\Pi_1\Pi_{12}\Pi_{13}\Pi_{14}$,
\[ \abs{1-\alpha_{12}-\alpha_{13}-\alpha_{14}}\leq 1.\]
Permuting the indices, we obtain the following $28$
inequalities in $12$ complex variables:
\begin{align}
\label{complextwo}
\abs{1-\alpha_{jk}}\leq 1,&
  \quad j,k\in\{1,2,3,4\} \text{ distinct},
\\
\label{complexthree}
\abs{1-\alpha_{jk_{1}}-\alpha_{jk_{2}}}\leq 1,&
  \quad j,k_{1},k_{2}\text{ distinct,}
\\
\label{complexfour}
\abs{1-\alpha_{jk_{1}}-\alpha_{jk_{2}}-\alpha_{jk_{3}}}\leq 1,&
  \quad j,k_{1},k_{2},k_{3}\text{ distinct.}
\end{align}
Let $\rho = e^{2\pi i/3} = (-1+i\sqrt{3})/2$,
a primitive cube root of unity, so $\bar\rho=\rho^2$.
We shall show:
\begin{lemma}
\label{lemmaC}
Complex numbers $\alpha_{jk}$ ($j,k\in\{1,2,3,4\}$ distinct) satisfy
\eqref{complexone}, \eqref{complextwo}, \eqref{complexthree},
and \eqref{complexfour} if and only if the equations
\begin{align*}
\alpha_{13}=\alpha_{31}=\alpha_{24}=\alpha_{42}& = 1,\\
\alpha_{12}=\alpha_{23}=\alpha_{34}=\alpha_{31}& = - \rho,\\
\alpha_{14}=\alpha_{43}=\alpha_{32}=\alpha_{21}& = - \bar\rho
\end{align*}
hold after permutation of the indices $1,2,3,4$.
\end{lemma}
This will suffice to prove Theorem~\ref{kelly-elkies},
because $\alpha_{13}=\alpha_{24}=1$ implies that
$\Pi_{13}$ and $\Pi_{24}$ are parallel, a contradiction.

We now prove the lemma.
\begin{proof}[Proof of Lemma~\ref{lemmaC}]
By \eqref{complexone}, we have as before
$\abs{\alpha_{jk}}+\abs{\alpha_{kj}}\geq 2$
with equality if and only if $\abs{\alpha_{jk}}=1$.
Summing over $j,k$, we obtain
\[
\sum_{\substack{j,k\\ j\neq k}}^3 \abs{\alpha_{jk}}\geq 12.
\]
We next prove that the sum is also bounded above by~$12$
by giving an upper bound on $\sum_{j\neq k}\abs{\alpha_{jk}}$
for each $j$.
\begin{lemma}
\label{sublemma}
Suppose $\alpha_1, \alpha_2,\alpha_3\in\C$ satisfy
\[
\Abs{1-\sum_{n\in S} \alpha_n}\leq 1,
\quad\text{for all $S\subseteq\{1,2,3\}$.}
\]
Then $\sum_{n=1}^3\abs{\alpha_n}\leq3$, with equality if and only if
$\{\alpha_1,\alpha_2,\alpha_3\}=\{1,-\rho,-\bar\rho\}$.
\end{lemma}
It will follow, from
\eqref{complextwo}, \eqref{complexthree}, and \eqref{complexfour}
and Lemma~\ref{sublemma},  that
\[
\sum_{\substack{j,k\\ j\neq k}}^3\abs{\alpha_{jk}}\leq 12.
\]
Thus equality holds, and again by Lemma~\ref{sublemma},
$\{\alpha_{jk}: k\neq j\}=\{1,-\rho,-\bar\rho\}$
for each $j=1,2,3,4$.
Keeping in mind that $\alpha_{jk}=\alpha_{kj}^{-1}$,
we obtain the $\alpha_{jk}$ values claimed in Lemma~\ref{lemmaC}.
\end{proof}

It remains to prove Lemma~\ref{sublemma}.
\begin{proof}[Proof of Lemma~\ref{sublemma}]

Since $\abs{1-\alpha_n}\leq 1$ for each $n$,
the $\alpha_n$ lie in the half plane $\operatorname{Re}(z) \geq 0$,
with $\operatorname{Re} (\alpha_n) > 0$ unless $\alpha_n = 0$.
Thus we may assume that the $\alpha_n$ are indexed so that
if none of them vanishes then
$\alpha_2$ lies in the angle $\myangle \alpha_1 0 \alpha_3$
determined by $\alpha_1$ and $\alpha_3$.
The given inequalities then imply that
\[
0,\
\alpha_1,\
\alpha_1+\alpha_2,\
\alpha_1+\alpha_2+\alpha_3,\
\alpha_2+\alpha_3,\
\alpha_3
\]
are the vertices (in this order) of a possibly degenerate hexagon
with perimeter $2\sum_{n=1}^3\abs{\alpha_n}$
contained in the disc $\abs{z-1}\leq 1$.
See Figure~\ref{fig2}.
\begin{figure}
\begin{center}
\begin{overpic}[scale=0.7]{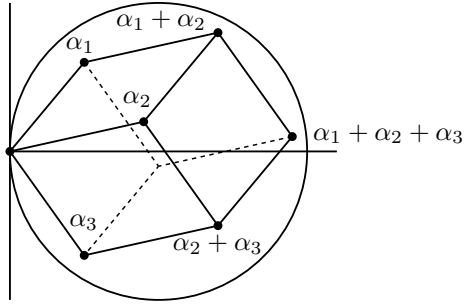}
\put(19,75){$\alpha_1$}
\put(36,59){$\alpha_2$}
\put(20,22){$\alpha_3$}
\put(34,82){$\alpha_1+\alpha_2$}
\put(93,48){$\alpha_1+\alpha_2+\alpha_3$}
\put(51,15){$\alpha_2+\alpha_3$}
\end{overpic}
\phantom{$\alpha_2+\alpha_3$}
\end{center}
\caption{Proof of Lemma~\ref{sublemma}}\label{fig2}
\end{figure}
However, it is easy to see that a hexagon contained
in a circle of radius $1$ has perimeter at most $6$,
with equality if and only if the hexagon is regular
and inscribed in the circle.  The lemma now clearly follows.
\end{proof}

This finishes the proof of Theorem~\ref{kelly-elkies}. \qed

\bigskip
Again we could have proved Theorem~\ref{kelly-elkies}
even if we allowed the plane at infinity so as to allow
parallel planes in~$\dualS$.
Then Lemma~\ref{lemmaC} would yield several tetrahedra
with the same volume as $\Pi_1\Pi_2\Pi_3\Pi_4$,
and iterating the argument produces an infinite (though locally finite)
set $\dualS_0 \subset \dualS$, contradicting the requirement
that $\dualS$ be finite, though $\dualS_0$ does satisfy
the dual SG incidence condition.
This configuration $\dualS_0$ is more easily described in affine
rather than barycentric coordinates for~$\C^3$.
Let $\cal E = \Z \oplus \Z\rho$,
the Eisenstein ring of algebraic integers in $\Q(\rho)$.
Then $\dualS_0$ can be identified with the set of lines in~$\C^3$
of the form
\[
x_1 - \rho^c x_2 = m,
\quad
x_2 - \rho^c x_3 = m,
\ \
\text{or}
\ \
x_3 - \rho^c x_1 = m
\]
with $c\in\{0,1,2\}$ and $m \in \cal E$.  For example,
if $x_1 - \rho^c x_2 = m$ and $x_2 - \rho^{c'} x_3 = m'$
then
$x_3 = \rho^{-c-c'} x_1 - \rho^{-c'}(m'+\rho^{-c}m)$,
and if $x_1 - x_2 = m$ and $x_1 - \rho x_2 = m'$ then
$x_1 - \bar\rho x_2 = -\rho m - \rho^2 m'$.
As in $\R^2$, we see that every configuration of planes in~$\C^3$
that satisfies the dual SG incidence condition
and has a tetrahedron of minimal volume is equivalent to $\dualS_0$
under some affine linear transformation of~$\C^3$.

\section{Proof of Theorem~\ref{main}}\label{quaternionsection}
It is sufficient to consider the four-dimensional case.
Dualizing the finite SG configuration spanning $\Proj^4(\Quat)$,
we obtain a finite collection $\CS$ of hyperplanes
with empty intersection, and with the property that
any two-dimensional flat in which two of the hyperplanes
in $\CS$ intersect, also contains a third hyperplane from $\CS$.
We call any five hyperplanes $\Pi_1, \Pi_2, \Pi_3, \Pi_4, \Pi_5$
with empty intersection a \emph{simplex}.
The \emph{vertices} of such a simplex are the five points
$P_i=\bigcap_{j\neq i}\Pi_j$, $i=1,\dots,5$.
Since all hyperplanes in $\CS$ have empty intersection,
$\CS$ contains at least one simplex.
Choose any hyperplane $\Pi_\infty$ at infinity
avoiding the vertices of some simplex in $\CS$.
We now define the \emph{measure} of a simplex $\Pi_1\dots\Pi_5$
to be the Study determinant of their vertices, i.e.,
\[ V(\Pi_1\Pi_2\Pi_3\Pi_4\Pi_5) := \sdet[P_1, P_2, P_3, P_4, P_5]\]
if all $P_p\notin\Pi_\infty$, $1\leq p\leq 5$,
and $V(\Pi_1\Pi_2\Pi_3\Pi_4\Pi_5) := \infty$ otherwise.
We now fix $\Pi_1\dots\Pi_5$ to be a simplex of minimum measure,
with vertices $P_1,\dots,P_5$.  (By the choice of $\Pi_\infty$,
there is at least one simplex of finite measure.
In fact we could have chosen any $\Pi_\infty$;
by the assumptions on $\CS$,
there will always be a simplex of finite measure.
Later we show that $\Pi_\infty$ can be chosen so that
all of those simplices have different measures,
and use such a choice to conclude the proof.)
We now change coordinates by a basis change in $\Quat^5$,
by letting each $P_i$ become the point associated with
the $1$-dimensional subspace of $\Quat^5$ generated by
the standard unit vector $e_p=[\delta_{p1},\dots,\delta_{p5}]^\Tr$,
and letting $\Pi_\infty$ become
the hyperplane associated with the $4$-dimensional subspace of $\Quat^5$
with equation $\sum_{p=1}^{5}x_p = 0$.  That is,
we now use barycentric coordinates with respect to $P_1,\dots,P_5$.
Because the Study determinant is multiplicative,
$\Pi_1\dots\Pi_5$ still has minimum measure.
Furthermore, $V(\Pi_1\Pi_2\Pi_3\Pi_4\Pi_5)=1$.

By the condition on $\CS$ we may choose, for each $p,q\in\{1,2,3,4,5\}$
with $p<q$, a hyperplane $\Pi_{pq}\neq\Pi_p,\Pi_q$
such that $\Pi_p\cap\Pi_q\subset\Pi_{pq}$.
Because $P_r\in\Pi_{pq}$ for all $r\neq p,q$
but $\Pi_{pq}\neq\Pi_p,\Pi_q$,
the equation of $\Pi_{pq}$ may be written as
\[ \alpha_{pq}x_p+x_q = 0\]
for some non-zero $\alpha_{pq}\in\Quat$.
We now set $\Pi_{qp}=\Pi_{pq}$,
and we can write its equation as $\alpha_{qp}x_q+x_p=0$ if we set
\begin{equation}\label{one}
\alpha_{qp}=\alpha_{pq}^{-1}.
\end{equation}

We now calculate the measures (all of which must be $\geq 1$)
of each of the following simplices:
\begin{enumerate}
\item $\Pi_1\Pi_{12}\Pi_3\Pi_4\Pi_5$,
\item $\Pi_1\Pi_{12}\Pi_{13}\Pi_4\Pi_5$,
\item $\Pi_1\Pi_{12}\Pi_{13}\Pi_{14}\Pi_5$,
\item $\Pi_1\Pi_{12}\Pi_{13}\Pi_{14}\Pi_{15}$.
\end{enumerate}

We consider the first case.
Let $\Pi'_p=\Pi_p$ for all $p\neq 2$, and $\Pi'_2=\Pi_{12}$.
Let $P_p'=\bigcap_{q\neq p}\Pi'_q$ for all $p$.
Then we calculate that
\[
P'_1=[(1-\alpha_{12})^{-1}, -\alpha_{12}(1-\alpha_{12})^{-1},0,0,0]^\Tr
\]
and $P'_p=e_p=P_p$ for each $p\geq 2$.
Then
\[
V(\Pi'_1\Pi'_2\Pi'_3\Pi'_4\Pi'_5) =
\sdet[P'_1,P'_2,P'_3,P'_4,P'_5] = \abs{(1-\alpha_{12})^{-1}}^2\geq 1,
\]
whence
\[ \abs{1-\alpha_{12}}\leq 1.\]
(In the case where $P'_1$ is at infinity and $V=\infty$ we have
$1-\alpha_{12}=0$, and the above inequality is satisfied trivially.)

Similarly, when doing the other cases, we obtain the following inequalities:
\begin{enumerate}
\item $\abs{1-\alpha_{12}}\leq 1$,
\item $\abs{1-\alpha_{12}-\alpha_{13}}\leq 1$,
\item $\abs{1-\alpha_{12}-\alpha_{13}-\alpha_{14}}\leq 1$,
\item $\abs{1-\alpha_{12}-\alpha_{13}-\alpha_{14}-\alpha_{15}}\leq 1$.
\end{enumerate}
By permuting indices we obtain $75$ inequalities in total:

\begin{align}
\abs{1-\alpha_{pq}} &\leq 1
 &\
 \text{($20$ inequalities)}
 \label{two}
\\
\abs{1-\alpha_{pq}-\alpha_{pr}}&\leq 1
 &\
 \text{($30$ inequalities)}
 \label{three}
\\
\abs{1-\alpha_{pq}-\alpha_{pr}-\alpha_{ps}}&\leq 1
 &\
 \text{($20$ inequalities)}
 \label{four}
 \\
\abs{1-\alpha_{pq}-\alpha_{pr}-\alpha_{ps}-\alpha_{pt}}&\leq 1
 &\
 \text{($5$ inequalities)}\label{five}
\end{align}
for any distinct $p,q,r,s,t\in\{1,2,3,4,5\}$.

We now prove the following:
\begin{lemma}\label{l1}
Assume $\alpha_{pq}\in\Quat$ ($p,q\in\{1,2,3,4,5\}$ distinct) satisfy
\eqref{one}, \eqref{two}, \eqref{three}, \eqref{four}, and \eqref{five}.
Then each $\alpha_{pq}=1/2 + \vv_{pq}$, where the vectors $\vv_{pq}$ satisfy
\begin{eqnarray}
\vv_{pq} &=& -\vv_{qp},\label{v1}\\
\ipr{\vv_{pq}}{\vv_{pq}} &=& \;\;\,3/4,\label{vnorm}\\
\ipr{\vv_{pq}}{\vv_{pr}} &=& -1/4,\label{vdot}
\end{eqnarray}
for all distinct $p,q,r\in\{1,2,3,4,5\}$,
and equality holds in all of
\eqref{two}, \eqref{three}, \eqref{four}, and \eqref{five}.
In other words, for each $p=1,\dots,5$ the set $\{\vv_{pq}:q\neq p\}$
comprises the vertices of a regular tetrahedron
inscribed in the sphere of radius $\sqrt{3}/2$ about the origin of~$\R^3$.
\end{lemma}

\begin{proof}
It follows from \eqref{one} and the AGM inequality that
\begin{equation}
\label{agm}
\abs{\alpha_{pq}}+\abs{\alpha_{qp}}\geq 2,
  \quad\text{with equality if and only if
  $\abs{\alpha_{pq}}=\abs{\alpha_{qp}}=1$.}
\end{equation}
Therefore, $\sum_{p,q:p\neq q}\abs{\alpha_{pq}}\geq 20$,
with equality if and only if each $\alpha_{pq}$ is a unit quaternion.

On the other hand we obtain an upper bound from the following lemma:

\begin{lemma}\label{l2}
Let $\beta_p\in\Quat, p=1,2,3,4$, satisfy
\begin{equation}
\label{star}
\Abs{1-\sum_{p\in A}\beta_p}\leq 1
  \quad \text{for all $A\subseteq\{1,2,3,4\}$.}
\end{equation}
Then $\sum_{p=1}^4\abs{\beta_p}\leq 4$,
with equality if and only if  each $\beta_p=1/2+\vv_p$,
where $\ipr{\vv_p}{\vv_p}=3/4$ and $\ipr{\vv_p}{\vv_q}=-1/4$
for all distinct $p,q$, in which case we also have equality
in each instance of \eqref{star}.
\end{lemma}

It will follow, by applying Lemma~\ref{l2}
to $\{\alpha_{pq} : q\neq p\}$
for each $p\in\{1,\dots,5\}$, that
$\sum_{p,q:p\neq q}\abs{\alpha_{pq}}\leq 20$.
Thus we have equality, and again by Lemma~\ref{l2} we obtain
\eqref{vnorm}, \eqref{vdot}, and equality in
\eqref{two}, \eqref{three}, \eqref{four}, \eqref{five}, and \eqref{agm}.
Since the $\alpha_{pq}$ are now unit quaternions,
we obtain \eqref{v1} from $\alpha_{qp}=\alpha_{pq}^{-1}=\alpha_{pq}^\ast$.
\end{proof}

\begin{proof}[Proof of Lemma~\ref{l2}]
In this lemma we use only the additive structure of the quaternions,
so we may consider the $\beta_p$ to be vectors $\vb_p\in\R^4$.
The given condition \eqref{star}
implies that the vertices of the parallelotope
\[
P := \Bigl\{\sum_{p=1}^4\lambda_p\vb_p:0\leq\lambda_p\leq 1\Bigr\},
\]
hence the whole $P$, is contained in the ball of radius~$1$
about $[1,0,0,0]^\Tr$.
The conclusion is equivalent to the statement that
the sum of the lengths of the four edges of~$P$\/
emanating from the vertex $\vo$ (the origin of~$\R^4$) is at most $4$,
with equality if and only if the $\vb_p$ form
an orthonormal basis of $\R^4$.  Therefore
the Lemma follows from the following slightly stronger statement:
\begin{quote}
  If $P$\/ is a parallelotope in $\R^4$ contained in a ball~$B$\/
  of unit radius, the sum of the lengths of the four generating vectors
  of~$P$\/ is at most $4$, with equality if and only if
  $P$\/ is a hypercube inscribed in the ball
  (and then necessarily of side length $1$).
\end{quote}
We now prove this statement.
Without loss of generality we may assume
$B$\/ is a unit ball about~$\vo$.  Let $\vc$ be the centroid of $P$.
Clearly $-P$, a parallelotope with centroid $-\vc$,
is also contained in~$B$.
Since $P$ is centrally symmetric, $-P$ is also a translate of $P$.
Since $B$\/ is convex, we may translate $P$\/ continuously
along a straight line to $-P$, with the translate staying inside $B$.
The centroids of these translates lie on the segment joining
$\vc$ and $-\vc$, hence one of these translates has centroid $\vo$.
Thus we have reduced the problem to parallelotopes with centroid~$\vo$.

Assume now that $P$\/ has centroid $\vo$.
Let $\va_p=\frac12\vb_p$ for each~$p$.
It is then enough to prove the following:
\begin{quote}
  If $\Norm{\sum_{p=1}^4\pm\va_p}\leq 1$ for all $2^4$ possible signs,
  then $\sum_{p=1}^4\norm{\va_p}\leq 2$,
  with equality if and only if
  the $\va_p$ are pairwise orthogonal and each of norm $1/2$.
\end{quote}
Indeed, if we sum
\begin{equation}
\label{ipr}
\Bigl\langle
  \sum_{p=1}^4\epsi_p\va_p,
  \sum_{q=1}^4\epsi_q\va_q
\Bigr\rangle\leq 1, \quad
(\epsi_1,\epsi_2,\epsi_3,\epsi_4)\in\{\pm 1\}^4,
\end{equation}
over all $16$ sign sequences,
we obtain $16\sum_{p=1}^4\norm{\va_p}^2\leq 16$,
whence the required $\sum_{p=1}^4\norm{\va_p}\leq 2$
follows by the Cauchy-Schwarz inequality.

Equality forces equality in Cauchy-Schwarz,
giving that all $\va_p$ have the same norm $\norm{\va_p}=1/2$,
and equality in \eqref{ipr}, giving
$\ipr{\va_p}{\va_q}=0$ for all distinct $p,q$.
\end{proof}

We have not quite completed the proof of Theorem~\ref{main}:
we have shown that the minimal simplex is not unique, but not
that one of its faces is parallel to another hyperplane in~$\dualS$.
But uniqueness suffices, because we could have chosen $\Pi_\infty$
so that no two simplices in~$\dualS$ have the same measure
unless that measure is infinite.  To see this,
note that there are finitely many pairs $\{\Sigma,\Sigma'\}$
of distinct simplices in~$\dualS$, and consider for each pair the set
$H(\Sigma,\Sigma')$ of hyperplanes $\Pi\subset\Proj^4(\Quat)$
such that $\Sigma$ and $\Sigma'$ have the same measure when
$\Pi$ is chosen as the hyperplane at infinity.
If we represent $\Pi$ using the five homogeneous coordinates
on the projective space dual to $\Proj^4(\Quat)$, and expand
each of these quaternions in its four real coordinates,
then the condition that $\Pi \in H(\Sigma,\Sigma')$
becomes a polynomial equation in these $5\times 4$ real variables.
Hence the union of our finitely many subsets $H(\Sigma,\Sigma')$
must have nonempty complement unless one of those polynomial equations
is satisfied identically.  But then $H(\Sigma,\Sigma')$
would consist of all the hyperplanes in $\Proj^4(\Quat)$.
This is absurd: since $\Sigma \neq \Sigma'$, some $P \in \Proj^4(\Quat)$
is a vertex of~$\Sigma$ but not of~$\Sigma'$;
and if $\Pi$ contains~$P$\/ but none of the other vertices
of $\Sigma$ and $\Sigma'$ then $\Sigma$ has infinite measure
while the measure of $\Sigma'$ is finite,
so $\Pi \notin H(\Sigma,\Sigma')$.
Therefore we may choose $\Pi_\infty$
outside the union of our sets $H(\Sigma,\Sigma')$.
The minimal simplex is then unique,
and applying Lemma~\ref{l1} to it yields a contradiction.
Our proof of Theorem~\ref{main} is now complete. \qed

\bigskip
Unlike the situation in $\R^2$ and $\C^3$,
in $\Quat^4$ there can be no infinite, periodic, locally finite
configuration of hyperplanes satisfying the dual SG condition.
We can show this by choosing a simplex
$\Pi_1 \Pi_2 \Pi_3 \Pi_4 \Pi_5$ of minimal measure
and checking that, for any arrangement of hyperplanes
around this simplex and the simplices
$\Pi_1 \Pi_{12} \Pi_3 \Pi_4 \Pi_5$ and $\Pi_1 \Pi_2 \Pi_{13} \Pi_4 \Pi_5$
of the same measure, the criteria of Lemma~\ref{l1}
must fail for at least one of these three simplices.
That computation also gives an alternative way
to deduce Theorem~\ref{main} from Lemma~\ref{l1}.

There is still a unique local configuration around
$\Pi_1 \Pi_2 \Pi_3 \Pi_4 \Pi_5$, analogous to those of
Lemmas \ref{trianglelemma} and~\ref{lemmaC},
that attains equality in each of the $75$ inequalities
\eqref{two}, \eqref{three}, \eqref{four}, and \eqref{five}.
It can be shown that \eqref{v1}, \eqref{vnorm}, and \eqref{vdot}
together imply that the $20$ vectors $\vv_{pq}$
are the vertices of a regular dodecahedron!
Starting from a regular dodecahedron inscribed in the sphere
of radius $\sqrt{3}/2$ about the origin of~$\R^3$,
we obtain $20$ unit quaternions by adding $1/2$ to each vertex;
there is then a unique way, up to permutations of $\{1,2,3,4,5\}$,
to set each of these quaternions equal to $\alpha_{pq}$
for some distinct $p,q \in \{1,2,3,4,5\}$ so that
$\alpha_{pq}+\alpha_{qp}=1$ and $|\alpha_{pq}-\alpha_{qr}|=1$
for all distinct $p,q,r \in \{1,2,3,4,5\}$.
Then equality holds in each of \eqref{two}, \eqref{three}, \eqref{four}, and \eqref{five}.
In $\R^2$ and $\C^3$, such a unique local configuration
led us to an infinite set of lines or planes satisfying
the dual SG incidence relation with a minimal triangle or tetrahedron.
What is the unique local configuration in~$\Quat^4$ hinting at?

\end{document}